\documentclass[12pt,draftcls,onecolumn]{IEEEtran}

\usepackage{amssymb,latexsym,epsfig}

\newcommand{\argmax}{\mathop{\mbox{\rm arg\,max}}}
\newcommand{\argmin}{\mathop{\mbox{\rm arg\,min}}}

\newtheorem{thm}{Theorem}[section]

\begin{document} 
\sloppy

\title{Sleeping Experts and Bandits Approach to
Constrained Markov Decision Processes}
 
\author{Hyeong Soo Chang\thanks{H.S. Chang is with the Department of Computer Science and Engineering at Sogang University, Seoul 121-742, Korea. (e-mail:hschang@sogang.ac.kr).}%
}

\maketitle

\begin{abstract}
This brief paper presents simple simulation-based algorithms for obtaining
an approximately optimal policy in a given finite set in large finite constrained Markov 
decision processes.
The algorithms are adapted from playing strategies for ``sleeping 
experts and bandits" problem and their computational complexities 
are independent of state and action space sizes if the given policy set is relatively small.
We establish convergence of their expected performances to the value of an optimal policy 
and convergence rates,
and also almost-sure convergence to an optimal policy with an exponential rate for 
the algorithm adapted within the context of sleeping experts.
\end{abstract}

\begin{keywords} 
constrained Markov decision processes, simulation, sleeping expert and bandit, 
learning algorithm
\end{keywords}

\section{Introduction}
Consider a discrete-time system with infinite horizon:
$
     x_{t+1}=f(x_t,a_t,w_t) \mbox{ for } t=0,1,2,...,
$ 
where
$x_t$ is the state at time $t$ 
-- ranging over a finite set $X$,
$a_t$ is the action at time $t$ 
-- to be chosen from a nonempty subset $A(x_t)$ of a given finite
set of available actions $A$ at time $t$, and $w_t$ is a random disturbance uniformly and independently selected from [0,1] at time $t$, representing the uncertainty in the system, and 
$f$ is a next-state function such that $f(x,a,w)\in X$ for $x\in X, a\in A(x)$, and $w\in [0,1]$.

Define a (stationary non-randomized Markovian) policy $\pi:X\rightarrow A$ with $\pi(x)\in A(x)$ for all $x\in X$ and
{\it value function} of $\pi$ given by
\begin{equation}
     V^{\pi}(x)=  E_{w_0,...,w_{\infty}} \biggl [ \sum_{t=0}^{\infty} \gamma^t R(x_t,\pi(x_t),w_t) \biggl | x_0 =x \biggr ],  x\in X, 
\label{valfun}
\end{equation}
with discount factor $\gamma\in (0,1)$ and one-period reward function 
$R$ such that $R(x,a,w) \in \mathcal{R}^+$ for $x\in X, a\in A(x)$, and $w\in [0,1]$ and
{\it constraint value function} of $\pi$ given by  
\begin{equation}
     J^{\pi}(x)=  E_{w_0,...,w_{\infty}} \biggl [ \sum_{t=0}^{\infty} \beta^t C(x_t,\pi(x_t),w_t) \biggl | x_0=x \biggr ],  x\in X,
\label{costvalfun}
\end{equation}
with discount factor $\beta\in (0,1)$ and one-period cost function 
$C$ such that $C(x,a,w) \in \mathcal{R}^+$ for $x\in X, a\in A(x)$, and $w\in [0,1]$. We let 
$R_{\max} = \sup_{x,a,w} R(x,a,w)$ and $C_{\max} = \sup_{x,a,w} C(x,a,w)$.

The function $f$, together with $X, A$, $R$, and $C$ comprise a constrained Markov decision process (CMDP)~\cite{altman}. For simplicity, we consider one constraint case. Extension
to multiple case is straightforward.

For a given $w=\{w_t\}$, we let $V^{\pi}(x,w)=\sum_{t=0}^{\infty} \gamma^t R(x_t,\pi(x_t),w_t)$ 
and $J^{\pi}(x,w)=\sum_{t=0}^{\infty} \beta^t C(x_t,\pi(x_t),w_t)$ with $x_0=x$. 
We assume throughout that any sample of $V^{\pi}(x,w)$ and $J^{\pi}(x,w)$ is bounded, respectively. Without loss of generality, we take the bound to be 1, 
i.e., for any $w$, $x$, and $\pi$,
$
    V^{\pi}(x,w) \in [0,1] \mbox{ and } J^{\pi}(x,w) \in [0,1].
$ (The generalization to an arbitrary bound can be done by appropriate scaling. Or by defining a transformation of $R$ into $R'$ such that $R'(x,a,w)=R(x,a,w)(1-\gamma)/R_{\max}$ and $C$ to $C'$ similarly, we can construct an ``equivalent" CMDP to the given CMDP which satisfies the assumption.) We also assume that an initial state $x_0$ is \emph{fixed} by some $x\in X$ and a nonempty finite policy set $\Pi$ is given.

A policy $\pi\in \Pi$ is called \emph{$\epsilon$-feasible} if $J^{\pi}(x) \leq K + \epsilon$ for given real constants $K>0$ and $\epsilon \geq 0$. 
We let $\epsilon$-feasible policy set $\Pi_{f}^{\epsilon} = \{\pi: \pi \in \Pi, J^{\pi}(x) \leq K + \epsilon \}$.
We then say that for $\epsilon\geq 0$, $\pi^{*}_{\epsilon}\in \Pi$ is an \emph{$\epsilon$-feasible optimal} policy 
if for some nonempty $\Delta$ such that $\Pi_{f}^{-\epsilon} \subseteq \Delta \subseteq \Pi_{f}^{\epsilon}$, $\pi^*_{\epsilon} \in \Delta$ and $\max_{\pi\in \Delta} V^{\pi}(x) = V^{\pi^*_{\epsilon}}(x)$.
The problem we consider is obtaining a $0$-feasible optimal policy (or estimating it with an $\epsilon$-feasible optimal policy) in $\Pi$, if such a policy exists.

The problem of obtaining a $0$-feasible optimal policy is known to be NP-hard if $\Pi$ contains all possible policies (in which case $|\Pi|=|A|^{|X|}$) and the problem size is characterized by the maximum of $|X|$ and $\max_{x\in X}|A(x)|$ and the number of constraints~\cite{feinberg}.
It seems that there exist only two exact iterative algorithms 
for this problem that exploit structural properties of CMDPs.
Chen and Feinberg~\cite{chen2} provided a value-iteration type algorithm
based on certain dynamic programming equations 
and Chang~\cite{changexact} presented a policy-iteration type algorithm based on a feasible-policy space characterization. 
Unfortunately, both require solving certain finite or infinite horizon MDP problems so that \emph{computational complexities 
depend on state and action space sizes.}
Note that linear programming used for finding a best \emph{randomized} policy cannot be applied here due to non-linearity and non-convexity of this problem (cf., P1 in~\cite[Theorem 3.1]{feinberg}). 

Even if there exists a body of works on simulation-based algorithms for solving unconstrained MDPs in order to break the curse of dimensionality (see, e.g.,~\cite{powell}~\cite{changbook} and the references therein), it seems that there has been no notable approach to CMDPs via simulation.
This paper is probably the first step toward developing such algorithms.
Because the algorithms proposed in this paper work with simulated sample-paths, computational complexities are independent of $|X|$ and $|A|$ as long as $|\Pi|$ is relatively small.

Our approach is simple and natural.
We generate a sequence of $\{\Pi_{f,n,H}, n=1,...,N\}$ where $\Pi_{f,n,H}$ is an estimate of $\Pi_f^0$, similar to the sample average approximation method~\cite{kleywegt}, by using simulation over a finite horizon $H$.
For each $\pi\in \Pi$, $J^{\pi}(x)$ is estimated with a sample mean and if the sample mean is less than or equal to $K$, $\pi$ is included in $\Pi_{f,n,H}$.
We then generate a sequence of policies $\{\pi(n), n=1,...,N\}$ from $\Pi_{f,n,H}$ at iteration $n$, where
$\pi(n)$ is an estimate of a 0-feasible optimal policy. 
The selection of $\pi(n)$ from $\Pi_{f,n,H}$ is based on the two playing strategies, called ``follow-the-awake-leader" (FTAL) and ``awake-upper-estimated-reward" (AUER),
for ``sleeping experts and bandits" problems~\cite{kleinberg}. 
A major difference between FTAL and AUER is that for FTAL, we simulate each policy in $\Pi_{f,n,H}$ to update the sample mean of each policy but for AUER, we simulate only selected policy $\pi(n)$ to update the sample mean of $\pi(n)$.
We view $\Pi_{f,n,H}$ as the set of currently awaken or non-sleeping experts/bandits in $\Pi$ and the sample value of the accumulated reward sum over the horizon $H$
as the sample reward of playing the expert/bandit $\pi$.
By proper adaptation of the results of the ``expected regret" defined over the sleeping experts and bandits model then,
we can establish convergence of the expected performance of our approach without the assumption that a 0-feasible optimal policy is unique.
We show that when $\Pi^0_{f} \neq \emptyset$, the expected performance $1/N\sum_{n=1}^N E[V^{\pi(n)}_H(x)]$ approaches the value of a 0-feasible optimal policy $\max_{\pi\in\Pi^0_f} V^{\pi}(x)$
as $N\rightarrow \infty$ and $H\rightarrow \infty$
with a rate of $O(1/N)$ (for $N\geq (\min_{\pi,\pi'\in \Pi} \{ V^{\pi}(x)-V^{\pi'}(x) : V^{\pi}(x)-V^{\pi'}(x) > 0\})^{-1}$) in the FTAL case and of $O(\ln N/N)$ in the AUER case for such $N$.
Here $V^{\pi}_H(x)=  E_{w_0,...,w_{H-1}} [ \sum_{t=0}^{H-1} \gamma^t R(x_t,\pi(x_t),w_t) | x_0 =x ]$ for $H<\infty$.
For the FTAL case, we further provide almost-sure convergence of $\pi(N)$ to a $0$-feasible optimal policy as $N$ and $H$ go to infinity with an exponential convergence rate at the expense of the assumption that value functions are all different among policies.

The works on the problem of finding the best solution from a finite set of solutions given
stochastic objective and constraint functions by simulation are relatively sparse 
(see~\cite{pasupathy} and the related references therein).
These works study allocating different (Monte-Carlo) simulation budgets to the solutions
to (approximately) maximize the probability of selecting the best solution from sample-mean estimates but provide explicit forms of such allocation only in an asymptotic limit, i.e., when the total 
number of samples approaches infinity. This is also typically given under the assumption that the 
best solution is unique and the distribution of samples are normal and in terms of the unknown 
true means and variances.
Even if heuristic iterative approximation procedures of such results are given, the convergences
of those are not known.
In our context, the best policy is not necessarily unique and the normality assumption is not
necessarily valid. 
Although Pasupathy et al.~\cite{pasupathy} consider general distribution case, the 
optimal allocation is only characterized by an optimization problem so that explicit forms of budget allocation
are difficult to obtain even in an asymptotic limit except for some special cases.
Without the uniqueness and the normality assumptions, Li et al.~\cite{li} consider a sequence of 
penalty cost functions to combine objective and constraint functions
with certain budget allocation strategy among the solutions but obtaining the sequence of the 
penalty cost functions is not straightforward
and their algorithm converges to a locally optimal solution when some restrictive assumptions 
are satisfied. 

Our setting also covers that in which explicit forms for $f$, $R$, and $C$ 
are not available, but they can be simulated. In this setting, another approach to consider is to
employ a stochastic-approximation based learning-algorithm as for unconstrained MDPs 
(see, e.g.,~\cite{djonin}~\cite{bhatnagar}). 
But this works when $\Pi$ is the set of all possible policies and the convergence speed is typically very slow and finite-time behaviours of such methods are not known.
Moreover, it's not immediate how to
adapt such approach when $\Pi$ is a subset of the set of all possible policies.

\section{Algorithm}
\label{sec:algo}

We first provide the pseudocode of the FTAL algorithm below.
It mainly consists of the \textbf{Feasible-Policy Set Estimation} step and 
the \textbf{Feasible Optimal Policy Estimation} step. The \textbf{Feasible-Policy Set Estimation} step obtains $\Pi_{f,n,H} = \{\pi : J^{\pi}_{n,H}(x) \leq K, \pi\in \Pi\}$ at iteration $n$. Here 
$J^{\pi}_{n,H}(x)$ is the sample mean obtained by $n$ independent samples of $J^{\pi}_H(x,w) = \sum_{t=0}^{H-1} \beta^t C(x_t,\pi(x_t),w_t)$ for $w=\{w_0,...,w_{H-1}\}$ and $H<\infty$.
We let $J^{\pi}_H(x) := E_{w} [ J^{\pi}_H(x,w) ]$.
The \textbf{Feasible-Policy Set Estimation} step selects $\pi(n)$ that achieves
$\max_{\pi\in \Pi_{f,n,H}} V^{\pi}_{\tau(\pi),H}(x)$ if $\Pi_{f,n,H} \neq \emptyset$ and $\tau(\pi)\neq 0$ for all $\pi \in \Pi_{f,n,H}$. (That is,
we ``follow the current best" among non-sleeping experts.)
Similarly, $V^{\pi}_{n,H}(x)$ is the sample mean obtained by $n$ independent samples of $V^{\pi}_H(x,w) = \sum_{t=0}^{H-1} \gamma^t R(x_t,\pi(x_t),w_t)$ for $w=\{w_0,...,w_{H-1}\}$, and $V^{\pi}_H(x) := E_{w} [ V^{\pi}_H(x,w) ]$.
The counter $\tau(\pi)$ keeps track of the number of times $\pi$ has been simulated to obtain a sample of $V^{\pi}_H(x,w)$. Whenever $\pi$ is included in $\Pi_{f,n,H}$ at some $n$, $\pi$ is simulated.
If there exists $\pi$ in $\Pi_{f,n,H}$ such that $\tau(\pi)=0$, $\pi(n)$ is set to be any such $\pi$.
If $\Pi_{f,n,H} = \emptyset$, $\pi(n)$ is set to be any $\pi \in \Pi$.

\noindent\textbf{Follow-The-Awake-Leader (FTAL)}
\begin{itemize} 
   \item[1.] \textbf{Initialization:} Select $N \geq 1$ and $H <\infty$. Set $J^{\pi}_{0,H}(x) = V^{\pi}_{0,H}(x) = 0$ and $\tau(\pi)=0$ for all $\pi\in \Pi$ and $n=1$.
   \item[2.] \textbf{Loop:} while ($n \leq N$)
       \begin{itemize}
          \item[2.1] \textbf{Feasible-Policy Set Estimation:} 
             For each $\pi\in \Pi$, obtain $J^{\pi}_H(x,w)$ by generating $w=\{w_0,...,w_{H-1}\}$ and set
\[
   J^{\pi}_{n,H}(x) = \frac{n-1}{n} J^{\pi}_{n-1,H}(x) + \frac{J^{\pi}_H(x,w)}{n}.
\] 
             Obtain $\Pi_{f,n,H} = \{\pi : J^{\pi}_{n,H}(x) \leq K, \pi\in \Pi\}$.
            \item[2.2] \textbf{Feasible Optimal Policy Estimation:} \\
           \textbf{If} ($\Pi_{f,n,H} \neq \emptyset$) \textbf{Then}
              \begin{itemize}
                 \item[-] \textbf{If} $\exists \pi \in \Pi_{f,n,H}$ such that $\tau(\pi)=0$, \textbf{Then} $\pi(n) = \pi$
                 \item[-] \textbf{Else} $\pi(n) \in \argmax_{\pi\in \Pi_{f,n,H}} V^{\pi}_{\tau(\pi),H}(x)$.
                 \item[-] For each $\pi \in \Pi_{f,n,H}$, obtain $V^{\pi}_H(x,w)$ by generating $w=\{w_0,...,w_{H-1}\}$ 
                and set
                \[
                   V^{\pi}_{\tau(\pi)+1,H}(x) = \frac{\tau(\pi)}{\tau(\pi)+1} V^{\pi}_{\tau(\pi),H}(x) + \frac{V^{\pi}_H(x,w)}{\tau(\pi)+1}
                \] and $\tau(\pi) \leftarrow \tau(\pi)+1$.
             \end{itemize}
	    \textbf{ElseIf} ($\Pi_{f,n,H} = \emptyset$) \textbf{Then} set $\pi(n)$ to be any policy in $\Pi$. 
          \item[2.3] $n \leftarrow n+1$
        \end{itemize}
\end{itemize}

As in FTAL, the AUER algorithm consists of the same two main steps. The \textbf{Feasible-Policy Set Estimation} step obtains $\Pi_{f,n,H}$ as in the FTAL case. 
Differently from the FTAL case, 
the \textbf{Feasible-Policy Set Estimation} step selects $\pi(n)$ at iteration $n$ which achieves
$\max_{\pi\in \Pi_{f,n,H}} ( V^{\pi}_{\tau(\pi),H}(x) + \sqrt{\frac{8\ln n}{\tau(\pi)}} )$ if $\Pi_{f,n,H} \neq \emptyset$ and $\tau(\pi)\neq 0$ for all $\pi \in \Pi_{f,n,H}$. 
(The term $\sqrt{\frac{8\ln n}{\tau(\pi)}}$ plays the role of estimating ``upper confidence bound" or ``upper estimated reward"~\cite{kleinberg}. We choose the bandit with the current highest upper estimated reward.) 
Then, only $\pi(n)$ is simulated and the sample mean of $\pi(n)$ is updated. The pseudocode of the AUER algorithm is given below.

\noindent\textbf{Awake-Upper-Estimated-Reward (AUER)}
\begin{itemize} 
\item[1.] \textbf{Initialization:} Same as FTAL
\item[2.] \textbf{Loop:} while ($n \leq N$)
\begin{itemize}
\item[2.1] \textbf{Feasible-Policy Set Estimation:} Same as FTAL
\item[2.2] \textbf{Feasible Optimal Policy Estimation:} \\
\textbf{If} ($\Pi_{f,n,H} \neq \emptyset$) \textbf{Then}
   \begin{itemize}
       \item[-] \textbf{If} $\exists \pi \in \Pi_{f,n,H}$ such that $\tau(\pi)=0$, \textbf{Then} $\pi(n) = \pi$.
       \item[-] \textbf{Else} $\pi(n) \in \argmax_{\pi\in \Pi_{f,n,H}} \biggl ( V^{\pi}_{\tau(\pi),H}(x) + \sqrt{\frac{8\ln n}{\tau(\pi)}} \biggr )$.
       \item[-] Obtain $V^{\pi(n)}_H(x,w)$ by generating $w=\{w_0,...,w_{H-1}\}$ and set
                \[
                   V^{\pi(n)}_{\tau(\pi(n))+1,H}(x) = \frac{\tau(\pi(n))}{\tau(\pi(n))+1} V^{\pi(n)}_{\tau(\pi(n)),H}(x) + \frac{V^{\pi(n)}_H(x,w)}{\tau(\pi(n))+1}.
                \] 
       \item[-] $\tau(\pi(n)) \leftarrow \tau(\pi(n))+1$
    \end{itemize}
\textbf{ElseIf} ($\Pi_{f,n,H} = \emptyset$) \textbf{Then} set $\pi(n)$ to be any policy in $\Pi$. 
\item[2.3] $n \leftarrow n+1$
\end{itemize}
\end{itemize}

\section{Convergence Analysis}
\label{sec.conv1}

We start with the convergence result of $\{\Pi_{f,n,H}\}$.
The following theorem establishes that as $N\rightarrow \infty$, $\Pi_{f,N,H}$ approaches
$\frac{\beta^H C_{\max}}{1-\beta}$-feasible policy set with the rate of $O(c(\epsilon,H)^N)$
for some constant $c(\epsilon,H) \in (0,1)$.
That is, $\Pi_{f,N,H}$ is arbitrarily
close to $\{\pi : J^{\pi}_H(x) \leq K, \pi\in \Pi \}$ as $N\rightarrow \infty$. 
By letting then $H\rightarrow \infty$ and $N\rightarrow \infty$, we can see that
it goes to the true feasible policy set $\Pi_f^0$ with an exponential convergence rate $O(\max\{ \beta^H, c(\epsilon,H)^N \})$.
We use the $O$-notation to mean that $f(x)=O(g(x))$ if there exist
real constants $M$ and $k$ such that $|f(x)|\leq M |g(x)|$ for all $x>k$ for $f:\mathcal{R}\rightarrow \mathcal{R}$ and $g:\mathcal{R} \rightarrow \mathcal{R}$.

\begin{thm}
\label{thm:feaset}
Let $\alpha_H=\frac{\beta^H C_{\max}}{1-\beta}$. Then for any $\epsilon > \alpha_H$,
\[
\Pr\{ \Pi_f^{-\epsilon} \subseteq \Pi_{f,N,H} \subseteq \Pi_f^{\epsilon} \} \geq  1-2|\Pi| e^{-2 (\epsilon - \alpha_H)^2N}.
\]
\end{thm}
\proof
The following proof is partly based on the proof of Proposition 1 in~\cite{wang}.
The complement of the event $\{\forall \pi\in \Pi, (\pi \in \Pi_f^{-\epsilon}\rightarrow \pi \in \Pi_{f,N,H}) \wedge (\pi\in \Pi_{f,N,H} \rightarrow \pi\in \Pi_f^{\epsilon})\}$ is
$\{\exists \pi\in \Pi, (\pi\in \Pi_f^{-\epsilon} \wedge \pi\notin \Pi_{f,N,H}) \vee  (\pi\in \Pi_{f,N,H} \wedge \pi \notin \Pi_f^{\epsilon}\}$.
This event is equal to
$\{\exists \pi\in \Pi, (J^{\pi}(x) \leq K-\epsilon \wedge J^{\pi}_{N,H} > K) \vee (J^{\pi}(x) > K +\epsilon \wedge J^{\pi}_{N,H}(x) \leq K)\}$, which is further equal to
$\{\exists \pi\in \Pi, (J^{\pi}_{N,H}(x) - J^{\pi}(x) > \epsilon) \vee (J^{\pi}_{N,H}(x) - J^{\pi}(x) < -\epsilon) \}$.
Therefore,
\begin{eqnarray*}
\lefteqn{\Pr\{ \Pi_f^{-\epsilon} \subseteq \Pi_{f,N,H} \subseteq \Pi_f^{\epsilon} \}}\\
& & \geq 1 - \Pr\{ \exists \pi\in \Pi, (J^{\pi}_{N,H}(x) - J^{\pi}(x) > \epsilon) \vee (J^{\pi}_{N,H}(x) - J^{\pi}(x) < -\epsilon) \} \\
& & \geq 1 - \Pr\{ \exists \pi\in \Pi, (J^{\pi}_{N,\infty}(x) - J^{\pi}(x) > \epsilon - \alpha_H) \vee (J^{\pi}_{N,\infty}(x) - J^{\pi}(x) < -\epsilon + \alpha_H) \} \\
& & \geq 1 - \sum_{\pi\in \Pi} \Pr\{ J^{\pi}_{N,\infty}(x) - J^{\pi}(x) > \epsilon - \alpha_H \} - \sum_{\pi\in \Pi} \Pr\{ J^{\pi}_{N,\infty}(x) - J^{\pi}(x) < -\epsilon + \alpha_H \},
\end{eqnarray*} where the second step follows from the fact that $-\alpha_H \leq J^{\pi}_{N,\infty}(x) - J^{\pi}_{N,H}(x)  \leq \alpha_H$.

Applying Hoeffding inequality~\cite{hoeff}, we finally have that
$
\Pr\{ \Pi_f^{-\epsilon} \subseteq \Pi_{f,N,H} \subseteq \Pi_f^{\epsilon} \} \geq  1-2|\Pi| e^{-2 (\epsilon - \alpha_H)^2N}.
$
\endproof

We remark that if $\Pi_{f}^{\epsilon}=\emptyset$ for $\epsilon > \alpha_H$, then $\Pi_{f,N,H}$ goes to the empty set as $N\rightarrow \infty$ so that when $\Pi_{f}^0=\emptyset$, $\Pi_{f,N,H}$ goes to the empty set as $N\rightarrow \infty$ and $H\rightarrow \infty$. That is, we can (approximately) identify the insolvability of the problem by these algorithms.
In what follows, we assume that $\Pi_f^0\neq \emptyset$.

\subsection{The FTAL algorithm performance}

We first establish almost-sure convergence of the FTAL algorithm. 
For this result, we need an assumption that $V^{\pi}(x) \neq V^{\pi'}(x)$ for all $\pi,\pi'\in \Pi$ for a technical reason. 

\begin{thm}
\label{thm:ftal}
Assume that $V^{\pi}(x) \neq V^{\pi'}(x)$ for all $\pi,\pi'\in \Pi$.
Let $r_H=\frac{\gamma^H R_{\max}}{1-\gamma}$ and $\alpha_H=\frac{\beta^H C_{\max}}{1-\beta}$
and $\Delta^{\pi}_N = \max_{\pi'\in \Pi_{f,N,H}}V^{\pi'}(x) - V^{\pi}(x)$ for nonempty $\Pi_{f,N,H}$ 
generated by FTAL.
Then for any $\epsilon > \alpha_H$,
\begin{eqnarray*}
\lefteqn{\Pr\{\max_{\pi\in \Pi_{f}^{-\epsilon}} V^{\pi}(x) \leq V^{\pi(N)}(x) \leq 
	\max_{\pi\in \Pi_{f}^{\epsilon}} V^{\pi}(x)\}} \\
	& & \hspace{3cm} \geq (1-2|\Pi| e^{-2 (\epsilon - \alpha_H)^2N}) \times (1- \sum_{\pi\in \Pi_{f,N,H} \setminus \{\pi^*_N\}}2e^{-2(\frac{\Delta^{\pi}_N}{2} - r_H)^2N}),
\end{eqnarray*} where $\pi^*_N \in \argmax_{\pi\in \Pi_{f,N,H}} V^{\pi}(x)$.
\end{thm}
\proof
From the assumption, $\pi^*_N$ is unique.
Then we have that
\begin{eqnarray*}
\lefteqn{\Pr\{\pi(N)\neq \pi^*_N \} \leq \sum_{\pi\in \Pi_{f,N,H} \setminus \{\pi^*_N\}} \Pr\{V^{\pi}_{N,H}(x) > V^{\pi^*_N}_{N,H}(x)\}}\\
& & \leq \sum_{\pi\in \Pi_{f,N,H} \setminus \{\pi^*_N\}} \biggl (\Pr\{ V^{\pi}_{N,H}(x) > V^{\pi}(x) + \frac{\Delta^{\pi}_N}{2} \} + \Pr\{ V^{\pi^*_N}_{N,H}(x) < V^{\pi^*_N}(x) - \frac{\Delta^{\pi}_N}{2} \} \biggr ) \\
& & \leq \sum_{\pi\in \Pi_{f,N,H} \setminus \{\pi^*_N\}} \biggl (\Pr\{ V^{\pi}_{N,\infty}(x) > V^{\pi}(x) + \frac{\Delta^{\pi}_N}{2} - r_H\} + \Pr\{ V^{\pi^*_N}_{N,\infty}(x) < V^{\pi^*_N}(x) - \frac{\Delta^{\pi}_N}{2} + r_H\} \biggr ) \\
& & \leq \sum_{\pi\in \Pi_{f,N,H} \setminus \{\pi^*_N\}}2e^{-2(\frac{\Delta^{\pi}_N}{2} - r_H)^2N} \mbox{ by Hoeffding inequality}.
\end{eqnarray*} The result follows then from
$\Pr\{\max_{\pi\in \Pi_{f}^{-\epsilon}} V^{\pi}(x) \leq V^{\pi(N)}(x) \leq \max_{\pi\in \Pi_{f}^{\epsilon}} V^{\pi}(x)\} \geq \Pr\{ \Pi_f^{-\epsilon} \subseteq \Pi_{f,N,H} \subseteq \Pi_f^{\epsilon} \} \times \Pr\{ \pi(N) = \pi^*_N \}.
$
\endproof

From the above theorem, we can see that $\pi(N)$ generated by the FTAL algorithm converges to a 0-feasible optimal policy as $N\rightarrow \infty$ and $H\rightarrow \infty$ almost surely if it is unique.

The theorem below establishes a finite-time bound on the expected performance of the FTAL algorithm without the assumption that the value functions of policies are different. 
The convergence of the expected performance of the FTAL algorithm follows then from this. 
Because the result is obtained by a direct application of the expected regret bound of the FTAL algorithm for sleeping experts~\cite[Theorem 6]{kleinberg}, a proof is omitted.

We construct a one-to-one mapping $I:\{1,2,...,|\Pi|\} \rightarrow \Pi$ such that $V^{I(i)} \geq V^{I(j)}$ for all $i,j\in \{1,...,\Pi\}$ with $i\leq j$.
For $y\geq 0$ and $i,j\in \{1,2,...,|\Pi|\}$, let $i_y(j) = \argmin\{i : i\leq j, \Delta_{I(i),I(j)}\leq y, i\in \{1,...,|\Pi|\} \}$ and 
$j_y(i) = \argmax\{j : j\geq i, \Delta_{I(i),I(j)}\leq y, j\in \{1,...,|\Pi|\} \}$, where $\Delta_{\pi,\pi'}=V^{\pi}_H(x) - V^{\pi'}_H(x)$ for $\pi,\pi'\in \Pi$. Note that we allow $\Delta_{\pi,\pi'} = 0$ for $\pi,\pi'\in \Pi$.
In what follows, the expectation is taken over the algorithm's random choices of $\{\pi(n)\}$
given a fixed sequence of $\{\Pi_{f,n,H}\}$.

\begin{thm}
\label{thm:exp1}
For every $\delta \geq 0$ and $\{\pi(n)\}$ generated by FTAL,
\begin{eqnarray*}
\lefteqn{0\leq \frac{1}{N} \sum_{n=1}^N \max_{\pi \in \Pi_{f,n,H}} V^{\pi}_H(x) -  \frac{1}{N} \sum_{n=1}^N E[V^{\pi(n)}_H(x)]}\\
& & \leq 2\delta +  
 \sum_{j=j_0(1)+1}^{|\Pi|} \frac{O(1)}{N \max\{\delta,\Delta_{I(i_0(j)-1),I(i_0(j))}\}} +
 \sum_{i=1}^{j_0(|\Pi|)-1} \frac{O(1)}{N \max\{\delta,\Delta_{I(j_0(i)),I(j_0(i)+1)}\}}
\end{eqnarray*}
 for any fixed sequence of $\{\Pi_{f,n,H}\}$ generated by FTAL.
\end{thm}

Note that by letting $N\rightarrow \infty$ and $H\rightarrow \infty$, $1/N\sum_{n=1}^N \max_{\pi \in \Pi_{f,n,H}} V^{\pi}_H(x)$ approaches arbitrarily close to $\max_{\pi \in \Pi_f^0} V^{\pi}(x)$ with an exponential rate $O(\max\{\gamma^H, \beta^H, c(\epsilon,H)^N \})$ by Theorem~\ref{thm:feaset}.
This implies that $\frac{1}{N} \sum_{n=1}^N E[V^{\pi(n)}_H(x)]$ approaches $\max_{\pi \in \Pi_f^0} V^{\pi}(x)$ as $N\rightarrow \infty$ and $H\rightarrow \infty$ with a rate of $O(1/N)$ for $N \geq (\min_{\pi,\pi'\in \Pi}\{\Delta_{\pi,\pi'}: \Delta_{\pi,\pi'}>0 \})^{-1}$ by setting $\delta = 1/N$.
In some sense, we can view the value of $\min_{\pi,\pi'\in \Pi}\{\Delta_{\pi,\pi'}: \Delta_{\pi,\pi'}>0 \})$ as the level of the difficulty of solving the problem. As it gets closer to zero, $N$ needs to get larger to obtain the rate.

\subsection{The AUER algorithm performance}

For the AUER algorithm, we are not be able to provide almost-sure convergence result as in Theorem~\ref{thm:ftal} for the FTAL algorithm. This is because it is difficult to establish
that
an upper bound on the probability of not choosing a 0-feasible optimal policy goes to zero
as $N\rightarrow \infty$ and $H\rightarrow \infty$ due to the term $\sqrt{\frac{8\ln n}{\tau(\pi)}}$. However, we can still provide the convergence of the expected performance of the AUER algorithm.
The following theorem establishes a finite-time bound on the expected performance of the AUER algorithm, again \emph{without the assumption that the value functions of policies are different}. 
As before, the result is from a direct application of the expected regret bound of the AUER algorithm for sleeping bandits~\cite[Theorem 12]{kleinberg}. 

\begin{thm}
\label{thm:exp2}
For every $\delta \geq 0$ and $\{\pi(n)\}$ generated by AUER,
\begin{eqnarray*}
\lefteqn{0\leq \frac{1}{N} \sum_{n=1}^N \max_{\pi \in \Pi_{f,n,H}} V^{\pi}_H(x) -  \frac{1}{N} \sum_{n=1}^N E[V^{\pi(n)}_H(x)]}\\
& & \leq 2\delta +  
 \sum_{j=j_0(1)+1}^{|\Pi|} \frac{O(\ln N)}{N \max\{\delta,\Delta_{I(i_0(j)-1),I(i_0(j))}\}} +
 \sum_{i=1}^{j_0(|\Pi|)-1} \frac{O(\ln N)}{N \max\{\delta,\Delta_{I(j_0(i)),I(j_0(i)+1)}\}}.
\end{eqnarray*}
 for any fixed sequence of $\{\Pi_{f,n,H}\}$ generated by AUER.
\end{thm}

From the above result, we see that $\frac{1}{N} \sum_{n=1}^N E[V^{\pi(n)}_H(x)]$ approaches $\max_{\pi \in \Pi_f^0} V^{\pi}(x)$ as $N\rightarrow \infty$ and $H\rightarrow \infty$ with a rate of $O(\ln N/N)$ for $N \geq (\min_{\pi,\pi'\in \Pi}\{\Delta_{\pi,\pi'}: \Delta_{\pi,\pi'}>0 \})^{-1}$ by setting $\delta = 1/N$.
Note that the rate of AUER is slower than FTAL's by a factor of $O(\ln N)$ at the expense of simulating only the selected policy at each iteration.

\section{Concluding Remarks}
\label{sec:conc}

Even if the discussions are made under the model of finite CMDPs, the proposed algorithms can be applied to CMDPs with infinite state and/or infinite action spaces as long as $\Pi$ is a finite set and each policy in $\Pi$ can be simulated. All of the results in the paper still hold in this case.

When we estimate feasible policy set in FTAL and AUER, we need to simulate all policies in $\Pi$. 
Developing a non-enumerative method for the feasible-policy set generation step is a good future work direction.

\end{document}